\documentclass[10pt,a4paper,twoside,leqno]{article}
\usepackage{amsfonts,amssymb}
\usepackage{eufrak}
\usepackage{amscd}

\newtheorem{defn}{Definition}

\vspace{0.5cm}

 4
\font\ebf=cmbx8
\font\erm=cmr8

\usepackage{graphicx}

\begin{document}

\thispagestyle{empty}

\noindent {\bf First observations on Prefab posets` Whitney numbers }

\vspace{0.7cm} {\it A. Krzysztof Kwa\'sniewski}

\vspace{0.2cm}

\noindent {\erm the Dissident - relegated by Bia\l ystok University authorities  }\\
\noindent {\erm from the Institute of Computer Science to Faculty of Physics}\\
\noindent {\erm ul. Lipowa 41,  15 424  Bia\l ystok, Poland}\\
\noindent {\erm e-mail: kwandr@gmail.com}\\

\vspace{0.2cm}

\noindent {\ebf Summary}
\noindent {\small We introduce a natural partial order   $\leq$ in
structurally natural  finite subsets of the cobweb prefabs sets
recently constructed by the present author. Whitney numbers of the
second kind of the corresponding subposet which constitute
Stirling-like numbers` triangular array - are then calculated and
the explicit formula for them  is provided. Next - in the second
construction - we endow the set sums of  prefabiants with such an
another partial order that their their Bell-like numbers include
Fibonacci triad sequences introduced recently by the present
author  in order to extend famous relation between binomial Newton
coefficients and Fibonacci numbers  onto the infinity of their
relatives among which there are also the Fibonacci triad sequences
and binomial-like coefficients (incidence coefficients included).
The first partial order is $F$-sequence independent while the
second partial order is $F$-sequence dependent where $F$ is the
so called admissible sequence determining cobweb poset by construction. An
$F$-determined cobweb poset`s Hasse diagram becomes Fibonacci tree
sheathed with specific cobweb if the sequence $F$ is chosen to be
just the Fibonacci sequence.

\vspace{0.2cm}
AMS Classification Numbers: 05C20, 11C08, 17B56 .
\vspace{0.1cm}

\vspace{0.1cm}

\noindent presented  (November $2006$) at the Gian-Carlo Rota
Polish Seminar\\
\noindent \emph{http://ii.uwb.edu.pl/akk/sem/sem\_rota.htm} 

\noindent published: Advances in Applied Clifford Algebras\\
Volume \textbf{18}, Number 1 / February, \textbf{2008},  57-73\\
ONLINE FIRST, Springer Link Date, Friday, August 10, 2007

\vspace{0.1cm}

\section{Introduction} The clue algebraic concept of combinatorics - the so called 
prefab (with associative and commutative composition) was introduced in
[1], see also [2,3]. The elements of prefabs are called since now on - prefabiants.
In [4] the present author had constructed a new broader class of prefab`s 
extending combinatorial structure based on the so called cobweb posets  (see Section 1.
[4] for the definition of a cobweb poset as well as a combinatorial interpretation of its
characteristic binomial-type coefficients - for example- fibonomial ones [5,6]).\\
Here  we introduce two natural partial orders: one  $\leq$ in
grading-natural subsets of cobweb`s prefabs sets [4] and in the second proposal
we endow the set sums of prefabiants with such another partial order that 
one may extend the Bell numbers to sequences of Bell-like numbers encompasing
among infinity of others the Fibonacci triad sequences introduced by the present
author in [7].

\section{Prefab based posets and their Whitney numbers.}
\vspace{1mm} 
Let the family $S$ of combinatorial objects
($prefabiants$) consists of all layers   $\langle\Phi_k
\rightarrow \Phi_n \rangle,\quad k<n,\quad k,n \in
N ={0,1,2,...}$ and an empty prefabiant $i$.

\noindent The set $\wp$ of prime objects consists of all
sub-posets $\langle\Phi_0 \rightarrow \Phi_m \rangle$  i.e. all
$P_m$`s $m \in N$ constitute from now on a family of prime $prefabiants$ which
we define after[4]in two steps. Namely accompanying the set $E$ of edges to the
set $V$ of vertices - one obtains the Hasse diagram where here down ${p,q,s}\in N $. 
(Convention: Edges stay for arrows directed - say - upwards - see examples below).
\begin{defn}
$$P=\langle V,E\rangle,\quad V=\bigcup_{0\leq p}\Phi_p ,\quad E
=\{\langle\langle j , p\rangle ,\langle q ,(p+1) \rangle
\rangle\}\bigcup\{\langle\langle 1 , 0\rangle ,\langle 1 ,1
\rangle \rangle\}, \quad where 1 \leq j \leq {p_F} , 1\leq q
\leq {(p+1)_F}.$$
\end{defn}.
The finite  cobweb sub-poset $P_m$ is then defined accodingly.
\begin{defn}
$P_m=\langle V_m,E_m\rangle,$ where $V_m = \bigcup_{0\leq s\leq m}\Phi_s$ and 
$E_m$  is defined as $E$  restricted to  $V_m$ by $1 \leq p \leq {m-1}$
is called the prime cobweb poset.
\end{defn}
Layer $$\langle\Phi_k \rightarrow \Phi_n \rangle$$ is considered
here to be the set of  all max-disjoint isomorphic copies (iso-copies)
of $P_m , m=n-k$ [4]. As a matter of illustration we quote after [4] examples
of cobweb posets` Hasse Diagrams [9] so that the layers become visualized.
\begin{center}
\includegraphics[width=100mm]{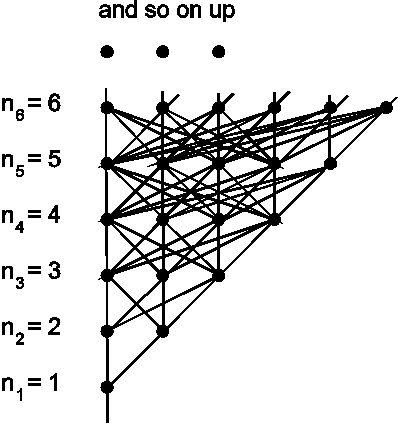}
\begin{center} 
{\small Fig.1. Display of Natural numbers` cobweb poset.}
\end{center}
\end{center}
\begin{center}
\includegraphics[width=100mm]{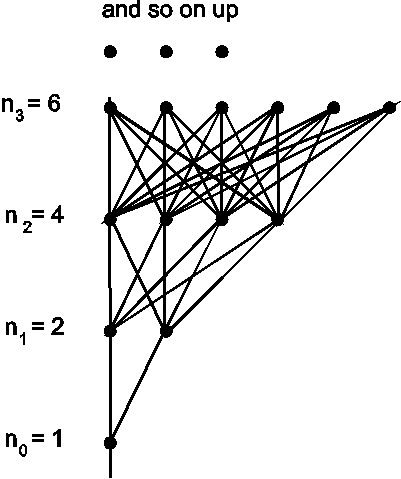}
\begin{center}{\small Fig.2. Display of Even Natural numbers $\cup \{1\}$-cobweb poset.}
\end{center}
\end{center}
\begin{center}
\includegraphics[width=100mm]{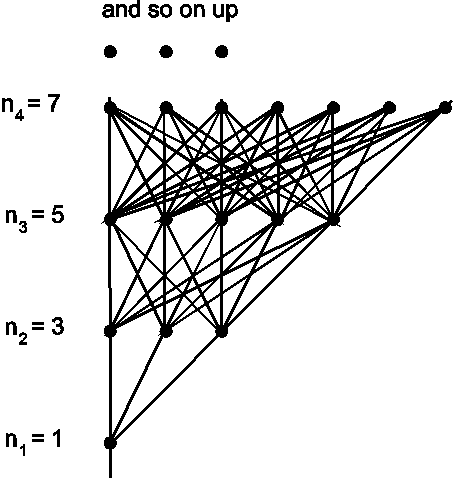}
\begin{center}{\small Fig3. Display of Odd natural numbers` cobweb poset.}
\end{center}
\end{center}
\begin{center}
\includegraphics[width=100mm]{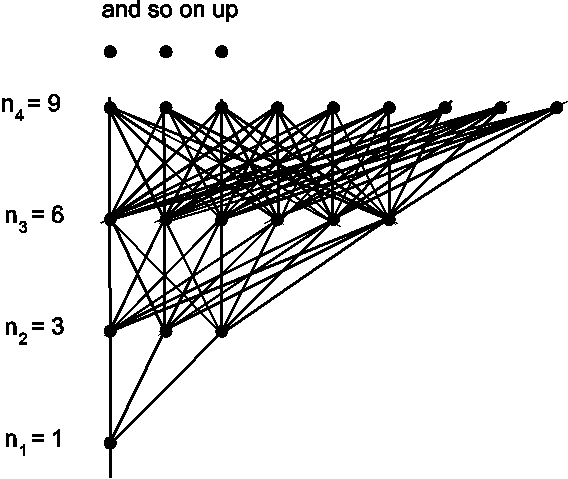}
\begin{center}{\small Fig.4. Display of divisible by 3 natural numbers $\cup \{1\}$ - cobweb poset.}
\end{center}
\end{center}
\begin{center}
\includegraphics[width=100mm]{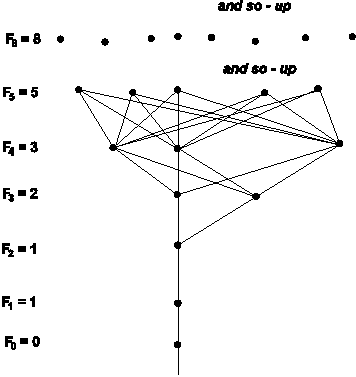}
\begin{center}{\small Fig.5. Display of Fibonacci numbers` cobweb poset.}
\end{center}
\end{center}

\section{Cobweb posets` combinatorial interpretation}
\noindent As seen above - for example  the $Fig.5$. displays the
rule of the construction of the  Fibonacci "cobweb" poset. It is
being visualized clearly while defining this [\textit{non-lattice!}] cobweb 
poset  $P$ with help of its incidence matrix [8]. The incidence $\zeta$
function  matrix representing uniquely just this cobweb poset $P$
has the staircase structure correspondent with $"cobwebed"$
Fibonacci Tree i.e. a Hasse diagram [9] of the particular partial
order relation under consideration.

\vspace{0.3cm}

$$ \left[\begin{array}{ccccccccccccccccc}
1 & 1 & 1 & 1 & 1 & 1 & 1 & 1 & 1 & 1 & 1 & 1 & 1 & 1 & 1 & 1 & \cdots\\
0 & 1 & 1 & 1 & 1 & 1 & 1 & 1 & 1 & 1 & 1 & 1 & 1 & 1 & 1 & 1 & \cdots\\
0 & 0 & 1 & 1 & 1 & 1 & 1 & 1 & 1 & 1 & 1 & 1 & 1 & 1 & 1 & 1 & \cdots\\
0 & 0 & 0 & 1 & 0 & 1 & 1 & 1 & 1 & 1 & 1 & 1 & 1 & 1 & 1 & 1 & \cdots\\
0 & 0 & 0 & 0 & 1 & 1 & 1 & 1 & 1 & 1 & 1 & 1 & 1 & 1 & 1 & 1 & \cdots\\
0 & 0 & 0 & 0 & 0 & 1 & 0 & 0 & 1 & 1 & 1 & 1 & 1 & 1 & 1 & 1 & \cdots\\
0 & 0 & 0 & 0 & 0 & 0 & 1 & 0 & 1 & 1 & 1 & 1 & 1 & 1 & 1 & 1 & \cdots\\
0 & 0 & 0 & 0 & 0 & 0 & 0 & 1 & 1 & 1 & 1 & 1 & 1 & 1 & 1 & 1 & \cdots\\
0 & 0 & 0 & 0 & 0 & 0 & 0 & 0 & 1 & 0 & 0 & 0 & 0 & 1 & 1 & 1 & \cdots\\
0 & 0 & 0 & 0 & 0 & 0 & 0 & 0 & 0 & 1 & 0 & 0 & 0 & 1 & 1 & 1 & \cdots\\
0 & 0 & 0 & 0 & 0 & 0 & 0 & 0 & 0 & 0 & 1 & 0 & 0 & 0 & 1 & 1 & \cdots\\
0 & 0 & 0 & 0 & 0 & 0 & 0 & 0 & 0 & 0 & 0 & 1 & 0 & 1 & 1 & 1 & \cdots\\
0 & 0 & 0 & 0 & 0 & 0 & 0 & 0 & 0 & 0 & 0 & 0 & 1 & 1 & 1 & 1 & \cdots\\
0 & 0 & 0 & 0 & 0 & 0 & 0 & 0 & 0 & 0 & 0 & 0 & 0 & 1 & 0 & 0 & \cdots\\
0 & 0 & 0 & 0 & 0 & 0 & 0 & 0 & 0 & 0 & 0 & 0 & 0 & 0 & 1 & 0 & \cdots\\
0 & 0 & 0 & 0 & 0 & 0 & 0 & 0 & 0 & 0 & 0 & 0 & 0 & 0 & 0 & 1 & \cdots\\
. & . & . & . & . & . & . & . & . & . & . & . & . & . & . & . & . \cdots\\
 \end{array}\right]$$

\vspace{1mm} \noindent \textbf{Figure 6.  The incidence matrix
$\zeta$ for the Fibonacci cobweb poset}

\vspace{2mm}

\noindent \textbf{Note} The knowledge of $\zeta$  matrix explicit
form enables one  to construct (count) via standard algorithms [8]
the M{\"{o}}bius matrix $\mu =\zeta^{-1} $ and other typical
elements of incidence algebra perfectly suitable for calculating
number of chains, of maximal chains etc. in finite sub-posets of
$P$. All elements of the corresponding incidence algebra are then
given by a matrix of the Fig.6 with $1$`s replaced by any reals (
or ring elements in more general cases).

\vspace{1mm}
We have a natural combinatorial object characterizig the cobweb posets
Hasse  directed graphs.

\vspace{1mm}

Namely - in general ([4-7], [10], [12], [13]) - given any sequence $\{F_n\}_{n\geq 0}$
of nonzero reals one may define  its corresponding  binomial-like $F-nomial$
coefficients in the spirit of \textbf{Ward`s Calculus of sequences}
[13]s as follows

\begin{defn}
$$
\left( \begin{array}{c} n\\k\end{array}
\right)_{F}=\frac{F_{n}!}{F_{k}!F_{n-k}!}\equiv
\frac{n_{F}^{\underline{k}}}{k_{F}!},\quad n_{F}\equiv F_{n}\neq
0, n\geq 0 $$ where we make an analogy driven identifications in
the spirit of  Ward`s Calculus of sequences  $(0_F\equiv0)$:
$$
n_{F}!\equiv n_{F}(n-1)_{F}(n-2)_{F}(n-3)_{F}\ldots 2_{F}1_{F};$$
$$0_{F}!=1;\quad n_{F}^{\underline{k}}=n_{F}(n-1)_{F}\ldots (n-k+1)_{F}. $$
\end{defn}
This is just the adaptation of the notation for the purpose
Fibonomial Calculus case (see references in [4-7], [10], 12]).
\vspace{2mm}

\noindent The crucial and elementary observation now is that the
cobweb poset combinatorial interpretation of $F$-binomial
coefficients [4-7,10,12,14,16] makes sense  \textbf{not for arbitrary} $F$
sequences as $F-nomial$ coefficients should be nonnegative
integers.

\vspace{2mm}

\begin{defn}
A sequence  $F = \{n_F\}_{n\geq 0}$ is called cobweb-admissible
iff
$$ \left( \begin{array}{c} n\\k\end{array}\right)_{F}\in N 
\quad for \quad k,n\in N .$$
\end{defn}

\vspace{2mm}

\noindent Right from the definition of $P$ via its
Hasse diagram  here now follow quite obvious and important
observations. They lead us to a combinatorial interpretation of
cobweb poset`s characteristic binomial-like coefficients (for
example - fibonomial ones [6,16]).  Here they are with the first
obvious observation  at the start.

\vspace{3mm}
14,
\noindent {\bf Observation 1}

{\it The number of maximal chains starting from The Root  (level
$0_F$) to reach any point at the $n-th$ level  with $n_F$ vertices
is equal to $n_{F}!$}.

\vspace{2mm}

\noindent {\bf Observation 2} $(k>0)$

{\it The number of all maximal chains in-between $(k+1)-th$ level $\Phi_{k+1}$
and the $n-th$ level $\Phi_n$ with $n_F$ vertices
is equal to $n_{F}^{\underline{m}}$, \quad where $m+k=n.$ } \\

\vspace{2mm}\noindent Indeed. Denote the number of ways to get
along maximal chains from  \textbf{any fixed point} (the leftist for example)
in $\Phi_k$ to $\Rightarrow  \Phi_n , n>k$ with the symbol\\
  $$[\Phi_k \rightarrow \Phi_n]$$
  then obviously we have :\\
           $$[\Phi_0 \rightarrow \Phi_n]= n_F!$$ and
$$[\Phi_0 \rightarrow \Phi_k]\times [\Phi_k\rightarrow \Phi_n]=
[\Phi_0 \rightarrow \Phi_n].$$

\vspace{2mm}

\noindent In order to formulate the combinatorial interpretation
of $F-sequence-nomial$ coefficients ({\it F-nomial} - in short)
let us consider all finite \textit{"max-disjoint"} sub-posets
rooted at the $k-th$ level at any \textbf{fixed} vertex $\langle
r,k \rangle, 1 \leq r \leq k_F $  and ending  at corresponding
number of vertices at the $n-th$ level ($n=k+m$) where the
\textit{"max-disjoint"} sub-posets are defined below.

\vspace{2mm}

\begin{defn}
Two families of maximal chains including two equipotent copies of $P_m$ 
are said to be max-disjoint if considered as sets of maximal chains 
they are disjoint i.e they have no maximal chain in common. (All $P_m$`s constitute
from now on a family of the so called  prime [4,10] $prefabiants$).  
An equipotent copy of $P_m$ [`\textbf{equip-copy}'] is  defined as such a family 
of maximal chains equinumerous with $P_m$ set of maximal chains that  the it constitues 
a sub-poset with one minimal element.
\end{defn}

\vspace{2mm}

\begin{defn}
We denote the number of all max-disjoint equipotent copies of
$P_m$  rooted at any vertex $\langle j,k \rangle , 1\leq j \leq
k_F $ of $k-th$ level  with the symbol
$$ \left( \begin{array}{c} n\\k\end{array}\right)_{F}.$$
\end{defn}
One uses the  customary convention:  $\left(
\begin{array}{c} 0\\0\end{array}\right)_{F}=1.$

\noindent Naturally- let us recall- the above definition makes sense 
\textbf{not for arbitrary} $F$ sequences as $F-nomial$ coefficients should be
nonnegative integers i.e. the sequence $F = \{n_F\}_{n\geq 0}$ must be cobweb-admissible.

\vspace{2mm}

\noindent \textbf{Problem 0.}
\textbf{The partition or tiling problem.}
Suppose now that  $F$  is a cobweb admissible sequence.
Let us introduce  

$$\sigma P_m = C_m[F; \sigma <F_1, F_2,...,F_m>].$$
the equipotent sub-poset obtained from $P_m$ with help of 
a permutation $\sigma$ of the sequence $<F_1, F_2,...,F_m>$. \\
Then 
              $$P_m = C_m[F; <F_1, F_2,...,F_m>].$$
Consider the layer  $\langle\Phi_k \rightarrow \Phi_n \rangle,\quad k<n,\quad
k,n \in N$.
Layer is considered here to be the set of  all max-disjoint
equipotent copies of $P_{n-k}$.  The question then arises , whether and under
which conditions  the layer may be partitioned with help of max-disjoint  blocks
of the form $\sigma P_m$.  At first - this main question answer is in affirmitave.
Some computer experiments done by student Maciej Dziemia\'nczuk [17] are encouraging 
However problems :  "how many?" or "find it all" are still \textbf{opened}.\\

\vspace{2mm}

Recall now that  the number of ways to reach an upper level from a
lower one along any of  maximal chains  i.e.  the number of all
maximal chains from the level
$\Phi_{k+1}$ to $ \Rightarrow  \Phi_n ,\quad n>k$ is equal to\\

  $$ [\Phi_k \rightarrow \Phi_n]= n_{F}^{\underline{m}}.$$

\noindent Naturally then we have

\begin{equation}
\left( \begin{array}{c} n\\k\end{array}\right)_{F} \times [\Phi_0
\rightarrow \Phi_m] = [\Phi_k \rightarrow \Phi_n]=
n_{F}^{\underline{m}}
\end{equation}
where  $[\Phi_0 \rightarrow \Phi_m]= m_F!$ counts the number of
maximal chains in any equip-copy of  $P_m$. With this in mind we see
that the following holds.

\vspace{3mm}

\noindent {\bf Observation 3} $(\textbf{n,k}\geq\textbf{0})$

{\it Let $n = k+m$. Let $F$ be any cobweb admissible sequence.
Then the number of max-disjoint  equip-copies i.e.  sub-posets
equipotent to $P_{m}$ , rooted at the same \textbf{ fixed} vertex of 
$k-th$ level and ending at the n-th level is equal to}

$$\frac{n_{F}^{\underline{m}}}{m_{F}!} =
\left( \begin{array}{c} n\\m\end{array}\right)_{F}$$
$$ = \left(\begin{array}{c} n\\k\end{array} \right)_{F}=
\frac{n_{F}^{\underline{k}}}{k_{F}!}. $$

\vspace{2mm} \noindent \textbf{Note} The above Observation 3 provides us
with the \textbf{new combinatorial interpretation} of  the  class
\textbf{of all classical $F-nomial$ coefficients} including
distinguished binomial or distinguished Gauss $q$- binomial ones
or Konvalina generalized binomial coefficients of the first and of
the second kind [11,12]- which include Stirling numbers too. The vast
family of Ward-like [13] admissible by
$\psi=\langle\frac{1}{n_F!}\rangle_{n\geq 0}$-extensions
$F$-sequences [12,14,16] includes also those desired here which shall
be called \textit{ "GCD-morphic"} sequences. This means that
$GCD[F_n,F_m] = F_{GCD[n,m]}$ where $GCD$ stays for Greatest
Common Divisor operator. The Fibonacci sequence is  a much
ontrivial [16,6] and guiding famous  example of GCD-morphic sequence.
Naturally  \textbf{incidence coefficients of any reduced incidence
algebra of full binomial type} [8] are GCD-morphic sequences
therefore they are now independently given a new cobweb
combinatorial interpretation via Observation 3. More on that - see
the next section where prefab combinatorial description is being
served. Before that - on the way - let us formulate  the following
problem (opened?).

\vspace{2mm}

\noindent \textbf{Problem 1} \textit{Find effective
characterizations of the cobweb admissible sequence i.e.  find all
examples. }

\noindent \textbf{Note on admissibility.} Observation 3 from  [16]
provides us with the\emph{ new combinatorial interpretation }of
the class \textbf{of all classical $F-nomial$ coefficients}
including distinguished binomial or distinguished Gauss $q$-
binomial ones or Konvalina generalized binomial coefficients of
the first and of the second kind [11]- which include Stirling
numbers too. This vast family of Ward-like [13] cobweb admissible
$F$-sequences - admissible at first by the so called
$\psi=\langle\frac{1}{n_F!}\rangle_{n\geq 0}$ umbral extensions[9]
- includes also those desired here which shall be called
\textit{ "GCD-morphic"} sequences.

\begin{defn}
The sequence of integers $F=\{n_F\}_{n\geq 0}$  is called the
GCD-morphic  sequence  if  $GCD[F_n,F_m] = F_{GCD[n,m]}$ where
$GCD$ stays for Greatest Common Divisor operator.
\end{defn}

\noindent Recall again : the Fibonacci sequence is  a much nontrivial [6] 
and guiding example of GCD-morphic sequence . Naturally incidence coefficients
of any reduced incidence algebra of full binomial type [8] are cobweb-admissible.
\noindent \textbf{Question:} which of these above are GCD-morphic
sequences?

\vspace{2mm}

\noindent In view of the \textbf{Note on admissibility} the
following problems are apparently interesting also on their own.

\vspace{2mm}

\noindent \textbf{Characterization Problem } \textit{Find
effective characterizations of the cobweb admissible sequence i.e.
find all examples.}

\noindent \textbf{GCD-morphism Problem } \textit{Find effective
characterizations i.e. find all examples. }

\vspace{2mm}

\section{Prefabs` Whitney numbers}

\noindent Consider then now  the partially ordered family $S$ of
these layers considered to be sets of  all max-disjoint isomorphic
copies (iso-copies) of  prime prefabiants $P_{m}= P_{n-k}$ as
displayed by Fig 1. - Fig.5. examples above. For any
$F$-sequence determining cobweb poset let us define in $S$
\textbf{the same} partial order relation as follows.

\begin{defn}
$$
\langle \Phi_k \rightarrow \Phi_n \rangle \leq
\langle\Phi_{k^\clubsuit} \rightarrow \Phi_{n^\clubsuit}\rangle
\quad \equiv\quad k \leq k{^\clubsuit} \quad\wedge\quad n\leq
n^{\clubsuit}.
$$
\end{defn}
For convenience reasons we shall also adopt and use the following
notation: $$ \langle\Phi_k \rightarrow \Phi_n \rangle = p_{k,n}.
$$
The interval $[p_{k,n}, p_{{k^\clubsuit},{n^\clubsuit}}]$  is of
course a subposet of $\langle\ S , \leq \rangle$. We shall
consider in what follows the subposet  $\langle\ P_{k,n} , \leq
 \rangle$ where
$$ P_{k,n}= [p_{o,o} ,p_{k,n}].$$

\noindent {\bf Observation 1.} The size $|P_{k,n}|$ of $P_{k,n}$ $
=|\{\langle l,m \rangle ,\quad 0 \leq l \leq k \quad\wedge \quad 0
\leq m \leq n \quad\wedge \quad k\leq n \}| = (n-k)(k + 1) +
\frac{k(k+1)}{2}.$

\vspace{2mm}

\noindent Proof: Obvious. Just draw the picture $\{\langle l,m
\rangle ,\quad 0 \leq l \leq k \quad\wedge \quad 0 \leq m \leq n
\quad\wedge \quad k\leq n \}$ of $P_{k,n}$` grid.

\vspace{1mm}

\noindent {\bf Observation 2.} The number of maximal chains in
$\langle\ P_{k,n} , \leq \rangle$ is equal to the number $d(k,n)$
of $0$ - dominated strings of binary i.e.  $0's$  and $1's$
sequences

$$d(k,n) = \frac{n+1-k}{n} \left( \begin{array}{c}
{k+n}\\n\end{array}\right). $$

\vspace{2mm}

\noindent Proof. The number we are looking for equals to the
number of minimal walk-paths in $[k\times n]$ Manhattan grid [15] -
paths restricted by the condition $ k \leq n $  i.e. it equals to
the number of $0$ - dominated strings of  $0's$  and  $1's$
sequences.

\vspace{2mm}

\noindent Recall that $( d(k,n) )$  infinite matrix`s  diagonal
elements are equal to the \textbf{Catalan} numbers $C(n)$

$$C(n) = \frac{1}{n} \left( \begin{array}{c}
{2n}\\n\end{array}\right). $$ as the Catalan numbers count the
number of  $ 0$ - dominated strings of   $0's$  and $1's$   with
equal number of $0's$  and $1's$ . Recall that  a $0$ - dominated
string of length  $n$  is such a string  that the first  $k$
digits of the string contain at least as many  $0's$ as $1's$ for
$k = 1, . . . , n $  i.e.  $0`s$ prevail  in appearance, dominate
$1`s$  from the left to the right end of the string.  $ 0$ -
dominated strings correspond bijectively  to minimal bottom - left
corner to the right upper corner paths in an integer grid $Z_\geq
\times Z_\geq$ rectangle part called 

Manhattan [15] with the restriction imposed on those minimal paths to obey
the "safety"  condition $k \leq n$ .

\vspace{2mm}

\noindent \textbf{Comment 1.} Observation 2. equips the poset
$\langle\ P_{k,n} , \leq \rangle$ with clear cut combinatorial
meaning.

\vspace{2mm}

\noindent The poset $\langle\ P_{k,n} , \leq \rangle$ is naturally
graded. $\langle\ P_{k,n} , \leq \rangle$ poset`s maximal chains
are all of equal size (Dedekind property) therefore the rang
function is defined.

\vspace{2mm}

\noindent {\bf Observation 3.} The rang $r(P_{k,n})$  of $P_{k,n}$
= number of elements in maximal chains $P_{k,n}$ minus one
$=k+n-1$. \noindent The rang $r(p_{l,m})$ of $\pi = p_{l,m} \in
P_{k,n}$ is defined accordingly: $r(p_{l,m})= l+m-1.$

\vspace{2mm}

\noindent Proof: obvious. Just draw the picture $\{\langle l,m
\rangle ,\quad 0 \leq l \leq k \quad\wedge \quad 0 \leq m \leq n
\quad\wedge \quad k\leq n \}$  of $P_{k,n}$` grid and note that
maximal means paths without at a slant edges.

\vspace{2mm}

\noindent Accordingly Whitney numbers $W_k(P_{l,m})$  of the
second kind are defined as follows (association: $ n
\leftrightarrow  \langle l,m \rangle $)

\vspace{2mm}

\begin{defn}
$$
W_k(P_{l,m})= \sum_{\pi\in P_{l,m}, r(\pi)=k} 1 \quad\equiv\quad
S(k,\langle l,m \rangle).
$$
\end{defn}

\vspace{2mm}

\noindent Here now and afterwords we identify  $W_k(P_{l,m})$ with
$S(k,\langle l,m \rangle)$ called and  viewed at  as Stirling -
like numbers of the second kind  of the naturally graded poset
$\langle\ P_{k,n} , \leq \rangle$ - note the association: $ n
\leftrightarrow  \langle l,m \rangle $.

\vspace{2mm}

\noindent \textbf{Right now challenge problems.  I.}

\vspace{2mm}

\noindent \textbf{I.} Let us define now Whitney numbers
$w_k(P_{l,m})$  of the first kind as follows (association: $ n
\leftrightarrow \langle l,m \rangle $.  Note the text-book
notation for M\"{o}bius function $\mu$)

\vspace{2mm}

\begin{defn}
$$
w_k(P_{l,m})= \sum_{\pi\in P_{l,m}, r(\pi)=k} \mu(0,\pi) \equiv
s(k,\langle l,m \rangle).
$$
\end{defn}

\vspace{0,2cm}

\vspace{2mm}

\noindent Here now and afterwards we identify  $w_k(P_{l,m})$ with
$s(k,\langle l,m \rangle)$  called and viewed at as Stirling -
like numbers of the first kind  of the poset $\langle\ P_{k,n} ,
\leq \rangle$ - note the association: $ n \leftrightarrow \langle
l,m \rangle $.

\noindent \textbf{Problem 1} Find an explicit expression for
$$
w_k(P_{l,m})\equiv s(k,\langle l,m \rangle) = ?
$$
and
$$
W_k(P_{l,m})\equiv S(k,\langle l,m \rangle) = ?
$$
Occasionally note that $S(k,\langle l,m \rangle)$ equals to the
number of the grid points counted at a slant (from the up-left to
the right-down) accordingly to the $l+m = k$ requirement.

\vspace{2mm}

\noindent \textbf{Problem 2} Find the recurrence relations for
$$
w_k(P_{l,m})\equiv s(k,\langle l,m \rangle \quad and\quad
W_k(P_{l,m})\equiv S(k,\langle l,m \rangle).
$$
We define now (note the association: $ n \leftrightarrow \langle
l,m \rangle $)  the corresponding Bell-like numbers
$$ B(\langle l,m \rangle)$$
of the naturally graded poset $\langle\ P_{k,n} , \leq \rangle$ as
follows.

\begin{defn}
$$
B(\langle l,m \rangle)= \sum_{k=0}^{l+m}S(k,\langle l,m \rangle).
$$
\end{defn}

\vspace{2mm}

\noindent {\bf Observation 4.}
$$ B(\langle l,m \rangle)=
|P_{l,m}| = \frac{k(k+1)}{2} + (n-k)(k+1).
$$
Proof: Just draw the picture $\{\langle l,m \rangle ,\quad 0 \leq
l \leq k \quad\wedge \quad 0 \leq m \leq n \quad\wedge \quad k\leq
n \}$ of $P_{k,n}$` grid and note that $S(k,\langle l,m \rangle)$
equals to the number of the grid points counted at a slant (from
the up-left to the right-down) accordingly to the $l+m = k$
requirement. Summing them up over all gives the size of $P_{k,n}$.

\vspace{2mm}
\noindent \textbf{Comment 2.} Observation 4. equips
the poset`s $\langle\ P_{k,n} , \leq \rangle$ Bell-like numbers
$B(\langle l,m \rangle)$ with clear cut combinatorial meaning.

\section{Set Sums of prefabiants` posets and their Whitney numbers.}

\vspace{1mm}

In this part we consider prefabiants` set sums with
an appropriate another partial order so as to arrive at Bell-like
numbers including Fibonacci triad sequences introduced recently by
the present author in [16] - see also [7,6].

\noindent Let $F$ be any   \textit{ "GCD-morphic"} sequence.
This means that $GCD[F_n,F_m] = F_{GCD[n,m]}$ where $GCD$ stays
for Greatest Common Divisor mapping. We define the
$F$-\textit{dependent} finite partial ordered set $P(n,F)$ as the
set of \textbf{prime} prefabiants $P_l$ given by the sum below.

\begin{defn}
$$ P(n,F) =\bigcup_{0\leq p}\langle \Phi_p \rightarrow \Phi_{n-p} \rangle = \bigcup_{0\leq l}P_{n-l}$$
\end{defn}
with the partial order relation defined  for $n-2l\leq 0$
according to
\begin{defn}
$$
P_l \leq P_{\hat l} \quad \equiv \quad l\leq \hat l,\quad P_{\hat
l}, P_l\in \langle \Phi_l \rightarrow \Phi_{n-l} \rangle.
$$
\end{defn}
Recall that \textbf{rang of} $P_l$ \textbf{is} $l$. Note that
$\langle \Phi_l \rightarrow \Phi_{n-l} \rangle = \emptyset$ for
$n-2l \leq 0$. The Whitney numbers of the second kind are
introduce accordingly.

\vspace{2mm}

\begin{defn}
$$
W_k(P_{n,F})= \sum_{\pi\in P_{n,F}, r(\pi)=k} \equiv S(n,k,F).
$$
\end{defn}
Right from the definitions above we infer that: (recall that
\textbf{rang of} $P_l$ \textbf{is} $l$.)

\vspace{2mm}

\noindent \textbf{Observation 5.}
$$
W_k(P_{n,F})= \sum_{\pi\in P_{n,F}, r(\pi)=k} \equiv S(k,n-k,F)=
\left( \begin{array}{c} {n-k}\\k\end{array}\right)_{F}.
$$
Here now and afterwords we identify $W_k(P_{n,F})=S(n,k,F)$ viewed
at and called as Stirling - like numbers of the second kind of the $P$ defined in [10]. 
$P$ by construction (see Figures above)
displays self-similarity property with respect to its prime
prefabiants sub- posets $P_n = P(n,F)$.

\vspace{2mm}

\noindent \textbf{Right now challenge problems. II.}

\noindent We repeat with obvious replacements of corresponding
symbols, names and definitions the same problems as in "Right now
challenge problems. I".

\noindent Here now consequently - for any $GCD$-morphic  sequence
$F$  (see: [10])  we define the corresponding Bell-like numbers $
B_n(F)$ of the poset $P(n,F)$ as follows.
\begin{defn}
$$B_n(F)= \sum_{k\geq0}S(n,k,F).$$
\end{defn}
Due to the investigation in [7,16] we have right now at our
disposal all corresponding results of [16,7] as the following
identification with  special case of $\langle
\alpha,\beta,\gamma\rangle$ - Fibonacci sequence $\langle
F_n^{[\alpha,\beta,\gamma]}\rangle_{n\geq 0}$ defined in [7]
holds.

\vspace{2mm}

{\bf Observation 6.}
$$B_n(F)\equiv F_{n+1}^{[\alpha=0,\beta=0,\gamma=0].}$$
\noindent Proof: See the Definition 2.2. from [7]. Compare also with the
special case of  formula (6) in [16].

\vspace{2mm}

\noindent\textbf{ Recurrence relations.} Recurrence relations for
$\langle \alpha,\beta,\gamma\rangle$ - Fibonacci sequences
$F_n^{[\alpha,\beta=,\gamma]}$ are to be found in [7] - formula
(9). Compare also with the special case formula (7) in [16].

\vspace{2mm}

\noindent\ \textbf{Closing-Opening Remark.} The study of further
properties of these Bell-like numbers as well as the study of
consequences of these identifications for the domain of the
widespread data types [7] and perhaps for eventual new dynamical
data types we leave for the possibly coming future. Examples of
special cases - a bunch of them - one finds in [7] containing [16]
as a special case. As seen from the identification Observation 6.
the special cases of $\langle \alpha,\beta,\gamma\rangle$ -
Fibonacci sequences $F_n^{[\alpha,\beta,\gamma]}$  gain
\textbf{additional} with respect to [16,7] combinatorial
interpretation in terms Bell-like numbers as sums over  rang $=k$
parts of the poset i.e. just sums of Whitney numbers of the poset
$P(n,F)$. This adjective \textit{"additional"} shines brightly
over Newton binomial connection constants between bases
$\langle(x-1)^k\rangle_{k\geq0}$ and $\langle x^n\rangle_{n\geq0}$
as these are Whitney numbers of the numbers from $[n]$ chain i.e.
Whitney numbers of the poset $\langle [n], \leq \rangle.$ For
other elementary "shining brightly"  examples see Joni , Rota and
Sagan excellent presentation in [18].

\section{ On applications of new cobweb posets` originated Whitney numbers}

\noindent Applications of new cobweb posets` originated    Whitney numbers  
such as extended Stirling or Bell numbers are expected to be of at least such
a significance in applications to linear algebra of formal series 
[linear algebra of generating functions  [19]] as Stirling and Bell numbers
or their $q$-extended correspondent already are in the so called  coherent physics
[20] ( see [20]  also for abundant references on the subject). 
Also straightforward applications of prefabs to coherent physics [20] are on line.
[Quantum coherent states physics  is of course a linear theory with its principle of states` superposition].

\vspace{2mm}

In order to say more on the subject of this section and give some examples let us remind
the equivalence of exponential structures by Stanley [21] with corresponding exponential prefabs [1]. 

\noindent In this context the let us indicate the crucial "Ward`ian - prefab`ian" 
example we owe to Gessel [22] with his $q$-analog of the exponential formula
as expressed by the Theorem 5.2  from [22].\\
We also recall that the $q$-analog of the Stirling numbers of the
second kind investigated by Morrison in Section 3  of  [23]
constitute {\bf the same} example of  Ward`ian  - prefab`ian
extension  as in the Bender - Goldman - Wagner  Ward - prefab
example. As noticed there by Morrison the  $(\gamma - e.g.f.)$
prefab exponential formula may equally well be derived from the
corresponding Stanley`s exponential formula in [21]. Let us then
now come over to these exponential structures of Stanley with an
expected impact on the current considerations  ( for definitions,
theorems etc. see [21]). In this connection we recall quoting
(notation from [21]) an important class of Stanley`s Stirling -
like numbers $\frac{S_{nk}}{M(n)}$ of the second and those of the
first kind Stanley`s Stirling - like numbers
$\frac{s_{nk}}{M(n)}$. Both kinds are characteristic immanent  for
counting of exponential structures (or equivalently -
corresponding exponential prefabs) and inheriting from there their
combinatorial meaning. This is due to the fact [21] that "with
each exponential structure is associated an "exponential formula"
and more generally a "convolution formula"  which is an analogue
of the well known exponential formula of enumerative
combinatorics" [21]. Consequently with each exponential structure
are associated Stirling-like , Bell-like numbers and Dobinski -
like formulas are expected also, of crucial impotance for "generalized coherent states` physics.

\vspace{2mm}
 As for the another examples let us consider in more detail exponential structures.

\vspace{2mm}
\noindent \textbf{Exponential structures.} Let $\{Q_n\}_{n\geq0}$
be any exponential structure and let  $\{M(n)\}_{n\geq0}$  be its
denominator sequence  i.e.  $M(n) =$ number of minimal elements of
$Q_n$. Let  $|Q_n|$  be the number of elements of the poset $Q_n$

$$ |Q_n| = \sum_{\pi\in Q_n} 1 .$$
Example: For $Q = \langle \Pi_n \rangle_{n\geq1}$ where $\Pi_n$ is
the partition lattice of  $[n]$ we have   $M(n)= 1$.

\noindent Define "Whitney-Stanley" number $S_{n,k}$ to be the
number of $\pi\in Q_n$ of degree  equal to $k\geq 1$ i.e.

$$ S_{n,k} = \sum_{\pi\in Q_n, |\pi|=k} 1 .$$
Define $S_{n,k}$ - generating characteristic polynomials (vide
exponential polynomials) in standard way

$$
W_n(x)= \sum_{\pi\in Q_n} x^{|\pi|}= \sum_{k=1}^{n} S_{n,k} x^k.
$$
Then  the exponential formula ($W_0(x)=1 = M(0)$) becomes

$$ \sum_{n=0}^{\infty}\frac{W_n(x) y^n}{M(n)n!} = exp\{xq^{-1}(y)\}, $$
where
$$q^{-1}(y) = \sum_{n=1}^{\infty}\frac{y^n}{M(n)n!}\equiv exp_{\psi} - 1 ,$$
with the obvious identification of $\psi$-extension choice here.
Hence the polynomial sequence $\langle p_n(x)= \frac
{W_n(x)}{M(n)}\rangle_{n\geq0}$ constitutes the sequence of
binomial polynomials  i.e. the basic sequence of the corresponding
delta operator $\hat Q = q(D)$. We observe then that
$$ p_n(x) = \sum_{k=0}^{n}\frac{S_{n,k}x^k}{M(n)}\equiv\sum_{k=0}^{n}[0,1,2,...,k;b_n]x^k
$$
are just exponential polynomials` sequence  for the equidistant
nodes case i.e. \textbf{Newton-Stirling} numbers of the second
kind $ S^\sim _{n,k}\equiv\frac{S_{n,k}}{M(n)}$. Both numbers and
the exponential sequence are being bi-univocally  determined by
the exponential structure $Q$. This is a special case of the one
considered in [20] and we have the - what we call- 
Newton-Stirling-Dobinski formula  (notation, history and details- see [20])
$$ p_n(x) = \frac {1}{\exp(x)}\sum_{k=0}^{\infty}\frac{b_n(k)
x^k}{k!}=
\sum_{k=0}^{n}[0,1,2,...,k;b_n]x^k,\quad\quad\quad\quad(N-S-Dob)$$
where $\langle b_n \rangle_{n\geq0}$ is defined by

$$b_n(x)=\sum_{k=0}^{n}S^\sim _{n,k}x^{\underline k}.$$
\textbf{Note} the \textbf{identification}
$b_n(x)=\frac{w_n(x)}{M(n)},$ where
$$
w_n(x)= -\sum_{\pi\in Q_n} \mu(\hat 0,\pi)\lambda^{|\pi|}.
$$
$\mu$ is M\"{o}bius function and $\hat 0$ is unique minimal
element adjoined to $Q_n$.

\noindent Corresponding Bell-like numbers [20] are then given by
$$ p_n(1) = \frac {1}{\exp(x)}\sum_{k=0}^{\infty}\frac{b_n(k)}{k!}=
\sum_{k=0}^{n}[0,1,2,...,k;b_n],\quad\quad\quad\quad(N-S-Bell).$$

\noindent Besides those above - in Stanley`s paper [21]  there are
implicitly present also inverse-dual "Whitney-Stanley" numbers
$s_{n,k}$ of the first kind i.e.

$$ s_{n,k} = -\sum_{\pi\in Q_n, |\pi|=k}\mu(\hat 0,\pi)  .$$
On this occasion and to the end of considerations on exponential
structures and Stirling like numbers let us make few remarks.
\noindent $q$-extension of exponential formula applied to
enumeration of permutations by inversions is to be find in
Gessel`s paper [22] (see there Theorem 5.2.) where among others he
naturally arrives at the $q$-Stirling numbers of the first kind
giving to them combinatorial interpretation. \noindent Recent
extensions of the exponential formula in the prefab language [1]
are to be find in [4]. Then \textbf{note}: exponential
structures, prefab exponential structures  (extended ones -
included) i.e. schemas where exponential formula holds-imply the
existence of Stirling like and Bell like numbers. As for the
Dobinski-like formulas one needs binomial or extended binomial
coefficients` convolution as it is the case with $\psi$-extensions
of umbral calculus in its operator form.

\vspace{0.2cm}

\textbf{Other Generalizations in brief.} We indicate here
\textsc{three} kinds of extensions of Stirling and Bell numbers -
including those which appear in coherent states` applications in
quantum optics on one side  or in the extended rook theory on the
other side. In the supplement for this brief account to follow on
this topics let us note that apart from applications to extended
coherent states` physics of quantum oscillators or strings  [6 -
11, 24, 25] and related Feymann diagrams` description [26] where
we face the spectacular and inevitable emergence of extended
Stirling and Bell numbers (consult  also [27]) there exists a good
deal of  work done on $discretization$ of space - time [28] and/or
Schrodinger equation using umbral methods [29] and GHW algebra
representations in particular (see: [28, 29] for references).\\

\vspace{2mm}

\noindent  \textbf{Acknowledgements}

\noindent Discussions with Participants of Gian-Carlo Rota Polish
Seminar  on all related topics \\
\emph{http://ii.uwb.edu.pl/akk/sem/sem\_rota.htm}  - are appreciated with
pleasure.

\begin
{thebibliography}{99}
\parskip 0pt

\bibitem{1}
E. Bender, J. Goldman   {\it Enumerative uses of generating
functions} , Indiana Univ. Math.J. {\bf 20} 1971), 753-765.

\bibitem{2}
D. Foata and M. Sch\"utzenberger, Th'eorie g'eometrique des
polynomes euleriens, (Lecture Notes in Math., No. 138).
Springer-Verlag, Berlin and New York, 1970.

\bibitem{3}
A. Nijenhuis and H. S. Wilf, Combinatorial Algorithms, 2nd ed.,
Academic Press, New York, 1978.

\bibitem{4}
A. K. Kwa\'sniewski, {\it Cobweb posets as noncommutative prefabs}
Adv. Stud. Contemp. Math.  vol. 14 (1), 2007 pp. 37-47. 

\bibitem{5}
A. K.  Kwa\'sniewski, {\it Information on combinatorial
interpretation of Fibonomial coefficients }   Bull. Soc. Sci.
Lett. Lodz Ser. Rech. Deform. 53, Ser. Rech.Deform. {\bf 42}
(2003), 39-41. ArXiv: math.CO/0402291   v1 18 Feb 2004

\bibitem{6}
A. K. Kwa\'sniewski, {\it The logarithmic Fib-binomial formula}
Advanced Stud. Contemp. Math. {\bf 9} No 1 (2004), 19-26. ArXiv:
math.CO/0406258 13 June 2004.

\bibitem{7}
A. K.  Kwa\'sniewski, {\it Fibonacci-triad sequences} Advan. Stud.
Contemp. Math. {\bf 9} (2) (2004),109-118.

\bibitem{8}
E. Spiegel, Ch. J. O`Donnell  {\it Incidence algebras}  Marcel
Dekker, Inc. Basel $1997$.

\bibitem{9}
$http://mathworld.wolfram.com/HasseDiagram.html$

\bibitem{10}
A. K. Kwa\'sniewski,  {\em Prefab posets` Whitney numbers }
Bull. Soc. Sci. Lett. Lodz, vol 55. 2005, pp. 17-25;  ArXiv:math.CO/0510027

\bibitem{11}
J. Konvalina , {\it A Unified Interpretation of the Binomial
Coefficients, the Stirling Numbers and the Gaussian Coefficients}
The American Mathematical Monthly {\bf 107} (2000), 901-910.

\bibitem{12}
A. K. Kwa\'sniewski {\it Main  theorems of extended finite
operator calculus} Integral Transforms and Special Functions, {\bf
14} No 6 (2003), 499-516.

\bibitem{13}
M. Ward: {\em A calculus of sequences}, Amer.J.Math. Vol.58,
(1936), 255-266.

\bibitem{14}
A. K. Kwa\'sniewski, {\it The logarithmic Fib-binomial formula}
Advanced Stud. Contemp. Math. {\bf 9} No 1 (2004), 19-26

\bibitem{15}
Z. Palka , A. Ruciñski {\it Lectures on Combinatorics}.I. WNT
Warsaw 1998 (\textit{in polish})

\bibitem{16}
A. K.  Kwa\'sniewski, {\it Fibonacci q-Gauss sequences} Advanced
Studies in Contemporary Mathematics {\bf 8} No 2 (2004), 121-124.
ArXive:  math.CO/0405591  31 May 2004.

\bibitem{17}
M. Dziemia\'nczuk {\it http:www.dejaview.cad.pl/brudnopisnaukowy.php}

\bibitem{18}
S.A. Joni ,G. C. Rota, B. Sagan {\it From sets to functions: three
elementary examples} Discrete Mathematics {\bf 37} (1981),
193-2002.

\bibitem{19}
H.S. Wilf, {\it Generatingfunctionology}    Boston: Academic Press,1990.

\bibitem{20}
A. K.  Kwa\'sniewski, {\it On umbral extensions of Stirling numbers and Dobinski-like formulas}
ASCM {\bf 12}(2006) no. 1, pp.73-100.

\bibitem{21}
R. Stanley, {\it Exponential structures}, Studies in Applied Math. {\bf 59} (1978), 73-82.

\bibitem{22}
I.M. Gessel {\it A q-analog of the exponential formula}  Discrete
Math. {\bf 40} (1982),  69-80

\bibitem{23}
Kent E. Morrison  {\it q-exponential families} The Electronic
Journal of Combinatorics {\bf 11} (2004) , No R36

\bibitem{24}
M. A. Mendez, P. Blasiak , K. A. Penson {\it Combinatorial
approach to generalized Bell and Stirling numbers and boson normal
ordering problem } arXiv : quant-ph/0505180  May   2005

\bibitem{25}
A.I. Solomon , P. Blasiak , G. Duchamp , A. Horzela , K.A. Penson
{\it  Combinatorial Physics, Normal Order and Model Feynman
Graphs} Proceedings of Symposium 'Symmetries in Science XIII',
Bregenz, Austria, 2003 arXiv: quant-ph/0310174 v1  29 Oct 2003

\bibitem{26}
Bender, C.M, Brody, D.C, and Meister, BK  {\it Quantum field
theory of partitions} Journal of Mathematical Physics, {\bf  40},
(1999),  3239-3245.

\bibitem{27}
Bender, CM, Brody, DC, and Meister, BK (2000) {\it Combinatorics
and field theory} Twistor Newsletter 45, 36-39.

\bibitem{28}
A Dimakis, F Müller-Hoissen and T Striker  {\it Umbral calculus,
discretization, and quantum mechanics on a lattice } J. Phys. A:
Math. Gen. {\bf 29}, (1996) 6861-6876

\bibitem{29}
D. Levi, P. Tempesta and P. Winternitz {\it Umbral Calculus,
Difference Equations and the Discrete Schroedinger Equation}
J.Math.Phys. {\bf 45} (2004) 4077-4105.

\end{thebibliography}



\end{document}